\documentclass[a4paper]{article}
\usepackage[british]{babel}

\usepackage{amsmath}
\usepackage{amssymb}
\usepackage{amsthm}

\usepackage{color}

\usepackage{hyperref}
\usepackage[nameinlink]{cleveref}

\newcommand{\mc}[1]{\mathcal{#1}}

\newcommand{\NN}{\mathbb{N}}
\newcommand{\ZZ}{\mathbb{Z}}
\newcommand{\QQ}{\mathbb{Q}}
\newcommand{\RR}{\mathbb{R}}
\newcommand{\vn}{{v_K}}
\newcommand{\bb}{\underline{b}}

\DeclareMathOperator{\res}{res}

\newtheorem{thrm}{Theorem}[section]
\newtheorem{prop}[thrm]{Proposition}
\newtheorem{lem}[thrm]{Lemma}
\newtheorem{corl}[thrm]{Corollary}
\newtheorem{fact}[thrm]{Fact}

\newcommand\blfootnote[1]{
	\begingroup
	\renewcommand\thefootnote{}\footnote{#1}
	\addtocounter{footnote}{-1}
	\endgroup
}

\theoremstyle{definition}

\newtheorem{defn}[thrm]{Definition}
\newtheorem{xmpl}[thrm]{Example}
\newtheorem{remk}[thrm]{Remark}

\title{Definable henselian valuations on dp-minimal real fields}
\author{Lothar Sebastian Krapp, Salma Kuhlmann, Lasse Vogel}
\date{}

\begin{document}
	\maketitle
	
	\begin{abstract}
		We give an explicit algebraic characterisation of all definable henselian valuations on a dp-minimal real field.
		Additionally we characterise all dp-minimal real fields that admit a definable henselian valuation with real closed residue field.
		We do so by first proving this for the more general setting of almost real closed fields.	\blfootnote{Math Subject Classification (2024): Primary 12J10 03C45; Secondary 03C60 03C40 13J30 16W60 12L12. Keywords: dp-minimal, almost real closed, definable valuation, henselian, real field.}		
	\end{abstract}
	
\section*{Introduction}
	
 In his seminal paper \cite[Conjecture 5.34~(c)]{Shelah14}, Shelah formulated a conjecture on algebraic properties of infinite strongly NIP fields (also termed strongly dependent). This conjecture has been adapted for the weaker condition NIP, a generalisation of strongly NIP, to the language of valuations. In this context the conjecture can be stated as:
	\emph{Every infinite NIP field is algebraically closed, real closed, or admits a non-trivial definable henselian valuation.}
	
	Various specialisations of this conjecture were considered in \cite{KrappKuhlmannLehericy21-2}, \cite{JahnkeSimonWalsberg17} and \cite{Johnson18}. In \cite{KrappKuhlmannLehericy21-2} the investigation restricts to strongly NIP, as in the original conjecture by Shelah. Strongly NIP imposes a boundary on the dp-rank (see \cite[Definition 4.12]{Simon15}), an important measure of complexity for NIP structures. In \cite{JahnkeSimonWalsberg17} and \cite{Johnson18} the query is further narrowed to fields which are dp-minimal (see \cite[Definition 4.27]{Simon15} for the definition). Dp-minimality is a further refinement of strongly NIP, limiting the dp-rank as much as possible.
	In \cite[Corollary 6.6]{JahnkeSimonWalsberg17} it is established that every dp-minimal \emph{ordered} field is real closed or admits a non-trivial definable henselian valuation.  Shelah's conjecture is verified for the dp-minimal case in \cite[Theorem 1.6]{Johnson18}, which establishes that every infinite dp-minimal field is real closed, algebraically closed or admits a non-trivial definable henselian valuation.\\
	
	In this note we are mainly interested in dp-minimal real fields, i.e.\ (formally) real fields $(K,+,\cdot)$ whose complete first-order theory is dp-minimal.
	From \cite{Johnson18} we immediately obtain that any dp-minimal \emph{real} field is \emph{either} real closed \emph{or} admits a non-trivial definable henselian valuation. Note that these two cases are exclusive, as a real closed field is o-minimal and therefore the only henselian valuation it defines is the trivial one (see \cite[Remark 5.1~(1)]{KrappKuhlmannLehericy21-2}).
	
	In \cite[Theorem 5.4]{KrappKuhlmannLehericy21-2}, Shelah's conjecture specialised to ordered fields is shown to be equivalent to the following: \emph{Every strongly NIP ordered field is almost real closed.}
	An adaptation to real fields can be found in \Cref{prop:dpMinRealFieldsClassification}, which shows that dp-minimal real fields are almost real closed.\\
	
	In the first part of this note (\Cref{section1}) we re-examine definable henselian valuations on almost real closed fields building on \cite{DelonFarre96}. In  \Cref{th:defblRCFresField} we obtain a criterion for the existence of a definable henselian valuation with real closed residue field on an almost real closed field, which in particular applies to dp-minimal fields due to \Cref{prop:dpMinRealFieldsClassification}.
	
	
	In the second part (\Cref{section2}) we focus our study entirely on dp-minimal real fields. We obtain the main result of this note (\Cref{thm:ClassificationOfDefinableHenselianValuations}), presenting a complete classification of definable henselian valuations on a dp-minimal real field in terms of the value group of its canonical henselian valuation.\\
	
	\noindent\textbf{Acknowledgements:} This work is part of the third author's doctoral research project `Dependent ordered structures', funded by \textsl{Evangelisches Studienwerk Villigst}. The first author was supported  by \textsl{Vector
	Stiftung} within the research project \textsl{Fundamentale Grenzen von Lernprozessen in künstlichen neuronalen Netzen, MINT-Innovationen 2022}.

\section{Definable valuations\newline on almost real closed fields}\label{section1}
	We start by gathering some preliminaries. Let $K$ be a (real) field and let $v$ be a valuation on $K$. We denote by
	\begin{itemize}
		\item $\mc O_v$ the valuation ring of $v$,
			
		\item $\mc M_v$ its maximal ideal,
			
		\item $v(K^\times)$ the value group, expressed as an additive group $(v(K^\times),+,<)$,
			
		\item  $Kv$ the residue field of $v$ and
			
		\item $\res_v\colon \mc O_v \longrightarrow Kv$ the residue map.
	\end{itemize}
	We say that a valuation $v$ is \textbf{definable} in $K$ if $\mc O_v$ is definable over $(K,+,\cdot)$ (with parameters). Usually we simply say that a set $A$ is definable in $K$ (respectively in $v(K^\times)$) if $A$ is definable over $(K,+,\cdot)$ (respectively $(v(K^\times),+,<)$).

	\begin{defn}\label{def:ARC}
		A field $K$ is called \textbf{almost real closed} if it admits a henselian valuation $v$ with real closed residue field $Kv$.
	\end{defn}
	
	Every almost real closed field $K$ is real, since due to the Baer--Krull Representation Theorem \cite[Theorem 2.2.5]{EnglerPrestel05} $K$ admits at least one ordering. By \cite[Theorem 5.2]{DittmannJahnkeKrappKuhlmann23}, if $(K,<)$ is an ordered almost real closed field, then every valuation definable in  $(K,+,\cdot,<)$ is henselian and already definable in $(K,+,\cdot)$. If on the other hand a valuation is definable in $(K,+,\cdot)$, then it is also definable in $(K,+,\cdot,<)$ for any ordering $<$ on $K$, so it follows that all definable valuations on an almost real closed field are henselian.
	By \cite[Proposition 2.1~(i)]{DelonFarre96} the set of henselian valuations on a real field $K$ is linearly ordered by inclusion of  the corresponding valuation rings (i.e.\ $v \leq w$ if and only if $\mc O_v \subseteq \mc O_w$). Hence any two henselian valuations are comparable. Furthermore the \textbf{canonical henselian valuation} $\vn$ (defined in \cite[page 106]{EnglerPrestel05}) is the minimal henselian valuation with respect to the ordering $\leq$. 
	\textbf{For a real field $K$ we fix the notation $G := \vn(K^\times)$}. By \cite[page 43~f.]{EnglerPrestel05}\ a convex subgroup $\Delta \subseteq G$ corresponds to a coarsening $v_\Delta\colon a \mapsto \vn(a) + \Delta$ of $\vn$ with value group $G/\Delta$. Note that $v_{\{0\}} = \vn$. From \cite[Corollary 4.1.4]{EnglerPrestel05} it follows that a coarsening $w$ of a henselian valuation $v$ is itself henselian. Hence $v_\Delta$ is a henselian valuation on $K$. Conversely for every henselian valuation $w$ on $K$ we obtain a convex subgroup $\Delta_w := \vn(\mc O_w \setminus \mc M_w) \subseteq G$ with $W = v_{\Delta_w}$. This yields a bijective correspondence of the henselian valuations on $K$ and the convex subgroups of $G$.
	
	\begin{defn}\label{def:SpecialValuations}
		Let $K$ be an almost real closed field. We denote the maximal convex divisible subgroup of $G$ as $G_0$. For a prime $p \in \NN$ we write $G_p$ for the maximal convex $p$-divisible subgroup of $G$. We set $v_0 = v_{G_0}$ and $v_p = v_{G_p}$.
	\end{defn}	
	
	Since divisibility implies $p$-divisibility $G_0 \subseteq G_p$ for every prime $p$.
	
	\begin{remk}\label{remark:v0isCoarsestWithRCres}
		By \cite[page 1123]{DelonFarre96} $v_0$ is the coarsest henselian valuation with real closed residue field on the almost real closed field $K$.
	\end{remk}

	We point out the following way to decompose the canonical henselian valuation.

	\begin{fact}\label{fact:valDecomp}
		\textup{\cite[page 44~f.]{EnglerPrestel05}}\ Let $K$ be real field and let $\Delta$ be a convex subgroup of $G$. Then $Kv_\Delta$ admits a henselian valuation $w$ such that $\vn$ is the composition of the valuations $v_\Delta$ and $w$. Furthermore $w({Kv_\Delta}^\times) = \Delta$. In particular if $K$ is almost real closed, then so is $Kv$.
	\end{fact}

	We now turn to establishing in \Cref{th:defblRCFresField} our characterisation of almost real closed fields, that admit a definable henselian valuation with real closed residue field. 

	\begin{prop}\label{prop:vpIsDefinable}
		Let $K$ be an almost real closed field and let $p$ be a prime. Consider the $\mc L_r$-formulae
		\begin{align*}
			\psi_p(x) &:=[ \neg\exists y\colon (y^p = x \vee y^p = -x) \wedge \exists z\colon z^p = 1+x ],\\
			\varphi_p(x) &:= \psi_p(x) \vee[\exists y\colon (y^p = x \vee y^p = -x) \wedge (\forall z\colon \psi_p(z) \rightarrow \psi_p(xz))] \vee x = 0.
		\end{align*}
		Then $\varphi_p(K) := \mc O_{v_p}$, i.e.\ the valuation $v_p$ is definable with defining formula $\varphi_p(x)$.
	\end{prop}
	\begin{proof}
		By applying \cite[Remark page 1126 and Proposition 2.6]{DelonFarre96} to $S = \{p\}$ in their notation, it suffices to show that $v_p = v_{\{p\}}$. We need to verfify
		\begin{itemize}
			\item[(1)] Hensel's lemma holds for polynomials $x^p - a$ for $a \in \mc O_{v_p}$,
			
			\item[(2)] $Kv_p = {Kv_p}^p \cup - {Kv_p}^p$ and
			
			\item[(3)] $v_p$ is the coarsest valuation on $K$ fulfilling (1) and (2).
		\end{itemize}
		(1) is clear since $v_p$ is henselian. By \cite[page 112]{DittmannJahnkeKrappKuhlmann23}, on any almost real closed field $K$ there exists a coarsest henselian valuation $w$ such that $Kw = {Kw}^p \cup -{Kw}^p$. By \cite[Lemma 5.1]{DittmannJahnkeKrappKuhlmann23} the corresponding convex subgroup $\Delta_w$ is the maximal $p$-divisible convex subgroup of $G$, i.e.\ $\Delta_w = G_p$ and since the correspondence is bijective it is $w = v_p$. This shows (2) and (3), completing the proof.
	\end{proof}

	Due to \Cref{prop:vpIsDefinable} all of the henselian valuations $v_p$ on the almost real closed field $K$ are definable. In the following we establish that a henselian valuation with real closed residue field is definable if and only if it is already one of the $v_p$.

	\begin{fact}\label{fact:defblConvSubgps}
		\textup{\cite[Corollary 4.3]{DelonFarre96}} Let $H$ be an ordered abelian group and $\Delta \not= \{0\}$ a definable convex subgroup of $H$. Then for some prime $p$ and $H_p$ the maximal $p$-divisible convex subgroup of $H$ it is $H_p \subseteq \Delta$.
	\end{fact}
	
	\begin{lem}\label{lem:canChooseNonNaturalDefblVal}
		Let $K$ be an almost real closed field and suppose that $\varphi(x;\bb)$ is a formula in the language $\{+,\cdot\}$ with parameters $\bb \in K^n$ defining $\vn$. Then there is an elementary extension $L \succ K$ such that $\varphi(L;\bb)$ is the valuation ring of a henselian valuation $w$ with real closed residue field $Lw$ and $w \not= v_L$.
	\end{lem}
	\begin{proof}
		Following the argument of \cite[Proof of Theorem 5.2, Case 2]{DittmannJahnkeKrappKuhlmann23}, there is an elementary extension $(L,+,\cdots,w) \succ (K,+,\cdot,\vn)$ such that $w$ is strictly coarser than $v_L$. It remains to note that $\varphi(L;\bb) = \mc O_{w}$. Since $\varphi(K;\bb) = \mc O_\vn$ it is $(K,+,\cdot,\vn) \models (\varphi(x;\bb) \leftrightarrow x \in \mc O_{\vn})$. Since $(L,+,\cdot,w)$ is an elementary extension it follows $ (L,+,\cdot,w) \models (\varphi(x;\bb) \leftrightarrow x \in \mc O_{w})$, so $\varphi(L;\bb) = \mc O_{w}$, which was to show.
	\end{proof}

	\begin{thrm}\label{th:defblRCFresField}
		Let $K$ be an almost real closed field. The following are equivalent:
		\begin{itemize}
			\item[(1)] $K$ admits a definable valuation $v$ with real closed residue field $Kv$.
			
			\item[(2)] There is a prime $p$ such that every $p$-divisible convex subgroup of $ G = \vn(K^\times)$ is already divisible.
			
			\item[(3)] $v_0$ is definable in $(K,+,\cdot)$.
		\end{itemize}
	\end{thrm}
	\begin{proof}\
		`(3) $\Rightarrow$ (1)': Follows immediately since by \Cref{remark:v0isCoarsestWithRCres} $Kv_0$ is real closed.\medskip
		
		\noindent `(2) $\Rightarrow$ (3)': By (2) there is some prime $p$ with $G_p = G_0$ and thus $v_p = v_0$. We now obtain (3) form \Cref{prop:vpIsDefinable}.\medskip
		
		\noindent `(1) $\Rightarrow$ (2)': Let $K$ be an almost real closed field and $\varphi(x;\bb)$ a formula such that $\varphi(K;\bb)$ is the valuation ring of a henselian valuation $v$ with real closed residue field. By \Cref{lem:canChooseNonNaturalDefblVal} consider an $L \succeq K$ such that $\mc O_w := \varphi(L;\bb)$ is not the valuation ring of the canonical henselian valuation. We denote the valuation with valuation ring $\mc O_w$ with $w$. Note that the theory of $K$ includes sentences witnessing that $\varphi(K;\bb)$ is the valuation ring of a henselian valuation with real closed residue field. It follows that $w$ is henselian and $Lw$ real closed.
		
		Let $G' := v_L(L^\times)$, $G'_0$ its maximal convex divisible subgroup and for a prime $p$ let $G'_p$ be the maximal convex $p$-divisible subgroup of $G'$.
		
		We want to show that for some prime $p$ every $p$-divisible convex subgroup of $G'$ is already divisible (i.e.\ $G'_p = G'_0$). In the structure $(L,+,\cdot,v_L)$ we can define $\Delta_w = v_L(\mc O_w \setminus \mc M_w)$, the convex subgroup of $G'$ corresponding to the coarsening $w$.
		
		By \cite[Corollary 5.25]{vandenDries14} the value group of a henselian valued field is stably embedded in the field. As a result it follows since we can define $\Delta_w \subseteq G'$ in the valued field $(L,+,\cdot,v_L)$, we can already define $\Delta_w$ in the ordered abelian group $(G',+,<)$. Now \Cref{fact:defblConvSubgps} yields that for some prime $p$ we have $G'_p \subseteq \Delta_w$. Note that every $p$-divisible convex subgroup of $G'$ is a subgroup of $G'_p$.
		
		On the other hand by \Cref{fact:valDecomp} $Lw$ admits a henselian valuation with value group $\Delta_w$. As $Lw$ was real closed, it follows from \cite[Theorem 4.3.7]{EnglerPrestel05} that $\Delta_w$ must be divisible. So every $p$-divisible convex subgroup of $G'$ is itself a subgroup of the divisible group $\Delta_w$ and therefore divisible.
		
		Now $w(L^\times) = G'/\Delta_w$ has no proper $p$-divisible convex subgroup. That means for all $x \in G'/\Delta_w$ there is $-|x| \leq y \leq |x|$ such that $y \not\in pG'/\Delta_w$. This is contained in the theory of $G'/\Delta_w$ as a sentence and by the Ax--Kochen--Ershov-Principle \cite[Theorem 5.1]{vandenDries14} it follows $G'/\Delta_w \equiv G/\Delta_v$. Hence $G/\Delta_v$ has no $p$-divisible convex subgroup which implies $G_p \subseteq \Delta_v$.
		
		But $\Delta_v$ is divisible by \cite[Theorem 4.3.7]{EnglerPrestel05} since $Kv$ is real closed and by \Cref{fact:valDecomp} there exists a henselian valuation with value group $\Delta_v$ on $Kv$. Hence $G_p$ and thus every $p$-divisible convex subgroup of $G$ is already divisible, which was to show.
	\end{proof}

	
	\begin{remk}\label{remk:v_0OnlyDefblWithRCFResField}
		Let $K$ be an almost real closed field. Further let $w$ be a definable henselian valuation with real closed residue field. Then $w = v_0$ follows from \cite[Proposition 5.9]{KrappKuhlmannLehericy21}. Also by \cite[Theorem 4.4]{DelonFarre96} this implies that $G_0$ is definable in $G$. Also note that whenever $G_0$ is non-trivial and definable in $(G,+,<)$, then it follows from \Cref{fact:defblConvSubgps}  that $G_0 = G_p$ for some prime $p$. In this case $v_0$ is definable since by \Cref{prop:vpIsDefinable} all $v_p$ are definable and recall that by \Cref{remark:v0isCoarsestWithRCres} the residue field $Kv_0$ is real closed.
	\end{remk}

	We conclude this section by giving an adaption of \Cref{th:defblRCFresField} on saturated almost real closed fields.
	
	\begin{corl}
		Let $K$ be an $\aleph_0$-saturated almost real closed field.\newline Then $K$ admits a definable valuation $v$ with real closed residue field if and only if the maximal convex divisible subgroup $G_0$ of $\vn(K^\times)$ is definable in $(\vn(K^\times),+,<)$.
	\end{corl}
	\begin{proof}
		The only thing that does not immediately follow from \Cref{remk:v_0OnlyDefblWithRCFResField} is: If $G_0 = \{0\}$, then $K$ admits a definable valuation $v$ with real closed residue field.
			
		
		So assume $G_0 = \{0\}$. Assume now for contradiction that furthermore $G_p \not= \{0\}$ for all primes $p$. Recall that $v_p$ is definable by \Cref{prop:vpIsDefinable}. The partial type
		\begin{align*}
			D(x) := & \{\forall y\colon [-|v_p(x)| \leq v_p(y) \leq |v_p(x)|\\
			& \Rightarrow (\exists z\colon pv_p(z) = v_p(y))] \mid p \textrm{ prime}\} \cup \{x \not= 0\}
		\end{align*}
		is consistent and thus realised in the saturated structure $(K,+,\cdot)$. Let $x_0$ be a realisation.
		As $x_0 \not= 0$, the convex subgroup of $G$ generated $v_K(x_0)$ is a non-trivial convex subgroup $\Delta$ of $G$. For every $g \in \Delta$ it is $g + G_p$ $p$-divisible in $G/G_p$, but since $G_p$ is $p$-divisible this implies that $g$ is $p$-divisible in $G$.
		Hence $\Delta$ is a non-trivial divisible convex subgroup of $G$ and therefore $\{0\} \subsetneq \Delta \subsetneq G_0 = \{0\}$. Contradiction.
	\end{proof}
	
	As the proof shows, this is possible since by assuming saturation the edge-case of the definable valuation being the canonical henselian valuation cannot occur. In \cite[Theorem 9.7]{Johnson23} the author obtains an analogue result. He does not specify to real fields and instead of almost real closedness requires a different property of the field $K$: dp-minimality. We will investigate this property for real fields in the following section.

\section{Classification of the definable valuations of a dp-minimal real field}\label{section2}

	The dp-rank \cite[Definition 4.12]{Simon15} is an important model theoretic measure of complexity. Roughly speaking, for a saturated model it is the supremum for the size of a family of mutually indiscernible sequences which can all be made discernible by adding a single constant to the language. A theory is called strongly NIP if in any saturated model there is no infinite family of mutually indiscernible sequences which can all be made discernible by adding a single constant to the language.
	
	Even more significant is the notion of dp-minimality \cite[Definition 4.27]{Simon15}, which means the dp-rank is 1, i.e\ for two mutually indiscernible sequences, adding a constant of the structure can not make both sequences discernible. In this section we now want to consider dp-minimal real fields.

	We first note that dp-minimal real fields are a special case of almost real closed fields. For dp-minimal fields a full characterisation is known, which we want to express specialised to real fields. For this we need the classification of dp-minimal ordered abelian groups.
	
	\begin{fact}\label{fact:dpMinOAGClassification}\textup{\cite[Proposition 5.1]{JahnkeSimonWalsberg17}}
		An ordered abelian group $H$ is dp-minimal if and only if $H/pH$ is finite for every prime $p$.
	\end{fact}
	
	With this we can now formulate the characterisation of dp-minimal real fields.
	
	\begin{prop}\label{prop:dpMinRealFieldsClassification}
		A real field $K$ is dp-minimal if and only if $K$ is an almost real closed field and $\vn(K^\times)/p\vn(K^\times)$ is finite for every prime $p$. 
	\end{prop}

	This proposition is analogous to \cite[Proposition 4.4]{KrappKuhlmannLehericy21-2}. It turns out we do not need to have a chosen ordering in the language.

	\begin{proof}
		By the last sentence of \cite[Theorem 1.2]{Johnson23} every dp-minimal field $K$ admits a henselian valuation $w$ whose residue field is finite, real closed or algebraically closed.  By \cite[Lemma 2.1]{KnebuschWright76} every henselian valuation on a real field is convex with respect to any ordering and thus has a real residue field. Then the residue field of $w$ is already real closed, hence $K$ is almost real closed.
		
		Now for a prime $p$ the valuation $v_p$ is definable, so the ordered abelian group $G/G_p$ is interpretable in $K$ and hence also dp-minimal. By \Cref{fact:dpMinOAGClassification} it follows that $(G/G_p)/p(G/G_p)$ is finite and $(G/G_p)/p(G/G_p) = G/pG$ as $G_p$ is $p$-divisible.
		
		For the other direction let $K$ be an almost real closed field such that the quotient $\vn(K^\times)/p\vn(K^\times)$ is finite for every prime $p$. By \Cref{fact:dpMinOAGClassification} $\vn(K^\times)$ is dp-minimal. Choose any ordering $<$ on $K$, then $(K,+,\cdot,<)$ meets the requirements in \cite[Proposition 4.4]{KrappKuhlmannLehericy21-2} and is therefore dp-minimal. Then so is the reduct $(K,+,\cdot)$.
	\end{proof}
	This implies in particular that if $K$ is a dp-minimal real field, then for any ordering $<$ on $K$ the ordered field $(K,<)$ is dp-minimal.\\
	
	Notably every dp-minimal real field is almost real closed, but not all almost real closed fields are dp-minimal. Consider the following example:
	
	\begin{xmpl}\label{xmpl:nonDpMinARCfield}
		We first define the following subgroup of $(\QQ,+,<)$:
		\[ B_0 := \left\{\frac{r}{s}\ \middle|\ r \in \ZZ, s \in \NN \textrm{ odd} \right\} \]
		Then $B_0$ is $p$-divisible for all primes $p \not= 2$. We can now consider $C_0 := B_0[\pi] \subsetneq \RR$ the additive group of the polynomials over $B_0$ evaluated at $\pi$. Then $\pi^i, \pi^j$ for $i \not= j \in \NN_0$  represent different elements of $C_0/2C_0$, hence this quotient is infinite and by \Cref{fact:dpMinOAGClassification} $C_0$ is not dp-minimal. Furthermore on the field of formal power series $\RR((C_0))$, the valuation given by the minimal exponent with non-zero coefficient is henselian, has residue field $\RR$ and value group $C_0$. By \Cref{prop:dpMinRealFieldsClassification} $\RR((C_0))$ is an almost real closed field that is not dp-minimal.
		
		Note that by \cite[Fact 4.6]{KrappKuhlmannLehericy21-2} the group $C_0$ is strongly NIP and by \cite[Lemma 4.8]{KrappKuhlmannLehericy21-2} so is the field $\RR((C_0))$.
	\end{xmpl}

	We now present the following examples to show that \Cref{th:defblRCFresField} is not void for dp-minimal real fields, i.e.\ there exist dp-minimal fields $K_1, K_2$, such that $K_1$ has the equivalent properties form \Cref{th:defblRCFresField} and $K_2$ does not.
	
		\begin{xmpl}\label{xmpl:defblARC}
			\begin{itemize}
				\item[(1)] Consider the group
				\[ C_1 := \ZZ \oplus \QQ \]
				with the lexicographic order. Then the only non-trivial proper convex subgroup is $\QQ$. Furthermore $C_1/pC_1 = \ZZ/p\ZZ$ for every prime $p$ and is therefore finite. So the group $C_1$ is dp-minimal by \Cref{fact:dpMinOAGClassification}. For any prime $p$ the maximal $p$-divisible subgroup of $C_1$ is $\{0\} \oplus \QQ$, hence it is already divisible. Therefore the field of formal power series $K_1 := \RR((C_1))$ has the equivalent properties from \Cref{th:defblRCFresField}. In particular $v_0 = v_p$ for all primes $p$ and can thus be defined with any of the formulae $\varphi_p(x)$ from \Cref{prop:vpIsDefinable}.

				\item[(2)] The group defined in \cite[Example 6.6~(ii)]{KrappKuhlmannLehericy21} was used as an example for a strongly NIP non-divisible ordered abelian group. Strongly NIP is a weaker property than dp-minimal: A structure is strongly NIP if the dp-rank of singletons is bounded by $\aleph_0$ \cite[Proposition 4.26]{Simon15}, opposed to dp-minmality asserting it to be $1$.
				We can show that this particular group $C_2$ is already dp-minimal:
				\[ \begin{array}{c}
					C_2 := \bigoplus_{k \in \NN} B_k\\
					B_n := \left\{ \frac{a}{p_{i_1}^{m_1} \cdot \ldots \cdot p_{i_k}^{m_k}} \middle| k \in \NN; i_1,\ldots,i_k \in \NN_0 \setminus \{n\}; a \in \ZZ; m_1,\ldots,m_k \in \NN_0 \right\}	
				\end{array}\]
				where $p_0 < p_1 < \ldots$ is an ordered list of all prime numbers in $\NN$.
				$C_2$ is $2$-divisible, so $C_2/2C_2 = \{0\}$. For any other prime $p_n$ the only component which is not $p_n$-divisible is $B_n$. As a result it is $C_2/p_nC_2 = B_n/p_nB_n = \ZZ/p_n\ZZ$ and therefore finite. As was pointed out in \cite{KrappKuhlmannLehericy21}, $C_2$ has no non-trivial divisible convex subgroup, but for every prime $p$ it has $p$-divisible convex subgroups. Hence an almost real closed field $K$ with $\vn(K^\times) = C_2$ fails to have the property (2) from \Cref{th:defblRCFresField}, and therefore also the others. The field $K_2 := \RR((C_2))$ is such an example and $K_2$ is dp-minimal by \Cref{prop:dpMinRealFieldsClassification}.
		\end{itemize}
	\end{xmpl}
	
	In \cite[Remark (4) after Theorem 4.4]{DelonFarre96} the authors state that under the condition that $G/pG$ is finite for all $p$ their results \cite[Theorem 4.4]{DelonFarre96} and \cite[Remark (1) after Theorem 4.1]{DelonFarre96} induce a full characterisation of real definable valuations. In our modern context this condition applies exactly to the dp-minimal real fields.
	
	We want to explicitly list all the definable henselian valuations of a dp-minimal real field. To show that we did not miss any definable henselian valuation we will need the following fact:
	
	\begin{fact}\label{fact:dfblConvSubgp_pRegularCondition}
	\textup{\cite[Theorem 4.1]{DelonFarre96}} Let $(H,<)$ be an ordered abelian group and let $\Delta \subseteq H$ be a definable convex subgroup. Then there is an $n \in \NN, n \not= 1$ such that for all convex subgroups $\Delta_1,\Delta_2 \subseteq H$ with $\Delta_1 \subsetneq \Delta \subsetneq \Delta_2$ the following holds: There are $g_1,g_2 \in \Delta_2/\Delta_1$ with $\exists h_1,\ldots,h_n \in \Delta_2/\Delta_1: g_1 \leq h_1 < \ldots < h_n \leq g_2$ and $\forall h: \neg(g_1 \leq nh \leq g_2)$ (i.e.\ $\Delta_2/\Delta_1$ is not $n$-regular; \cite[Appendix]{DelonFarre96}). Alternatively, there is a non-trivial convex subgroup $\Delta' \subsetneq \Delta_2/\Delta_1$ such that $(\Delta_2/\Delta_1)/\Delta'$ is not $n$-divisible. 
		
		We can always assume $n$ to be a prime by \cite[Theorem A~(vi)]{DelonFarre96}
	\end{fact}
	
	We can now give an explicit description of any definable henselian valuation on such a field:

	\begin{thrm}\label{thm:ClassificationOfDefinableHenselianValuations}
		Let $K$ be a dp-minimal real field and let $G := \vn(K^\times)$. Define $n_p \in \NN_0$ such that $|G/pG| = p^{n_p}$. There exists a surjective map
		\[\begin{array}{ccc}
			\{(p,i) \in \NN\times\NN_0 \mid p\textrm{ prime}, i \leq n_p\} & \rightarrow & \textrm{definable henselian valuations on K}\\
			(p,n) & \mapsto & v_{(p,n)}
		\end{array}\]
		where $v_{(p,n)}$ is given by composition of $\vn\colon K^\times \longrightarrow G$ with the projection $G \longrightarrow G/G_{(p,n)}$ with $G_{(p,n)}$ being the maximal convex subgroup $H \subseteq G$ such that $|H/pH| \leq p^n$.
	\end{thrm}
	\begin{proof}
		It suffices to show:
		\begin{itemize}
			\item[(1)] The $v_{(p,n)}$ are definable.
			
			\item[(2)] Every definable valuation is of the shape $v_{(p,n)}$.
		\end{itemize}
		(1): We first note that $G_{(p,0)} = G_p$ for every $p$ by its definition. As a result $v_{(p,0)} = v_p$ is definable by \Cref{prop:vpIsDefinable}. For $n \geq 1$ it is $G_{(p,0)} \subseteq G_{(p,n)}$, so $v_{(p,n)}$ is a coarsening of $v_{(p,0)}$. It is
		\[ \mc O_{v_{(p,n)}} = \mc O_{v_p} \cup \{x \mid v_p(x) \in G_{(p,n)}/G_{(p,0)} \}, \]
		so it suffices to show that $S_{(p,n)} := \{x \mid v_p(x) \in G_{(p,n)}/G_{(p,0)} \}$ is definable.
		
		To define this set choose parameters $x_1,\ldots,x_{p^n} \in K$ such that
		\[ \{ \vn(x_i) + p(G_{(p,n)}) \mid 1 \leq i \leq p^n \} = G_{(p,n)}/p(G_{(p,n)}), \]
		which is possible because $|G_{(p,n)}/p(G_{(p,n)})| \leq p^n$. Now consider the following formula:
		\[\small\begin{array}{rc}
			& \left[\varphi_p(x) \wedge \left(\forall y \not= 0\colon \left(\varphi_p(y) \wedge \varphi_p(\frac{x}{y})\right) \Rightarrow \displaystyle \bigvee_{i=1}^{p^n} \left(\exists z\colon \varphi_p(\frac{x_iy}{z^p}) \wedge \varphi_p(\frac{z^p}{x_iy})\right)  \right)\right]\smallskip\\
			\vee & \left[\neg\varphi_p(x) \wedge \left(\forall y \not= 0\colon \left(\neg\varphi_p(y) \wedge \varphi_p(\frac{y}{x})\right) \Rightarrow \displaystyle \bigvee_{i=1}^{p^n} \left(\exists z\colon \varphi_p(\frac{x_iy}{z^p}) \wedge \varphi_p(\frac{z^p}{x_iy})\right)  \right)\right]
		\end{array}\]
		We denote this formula by $\psi_{(p,n)}(x)$. An $x \in K^\times$ fulfils this formula if and only if for every $y$ with $-|v_p(x)| \leq v_p(y) \leq |v_p(x)|$ it is $v_p(y) + v_p(x_i) = v_p(yx_i)$ $p$-divisible for at least one of the $x_i$. Since $v_p(yx_i) = v_K(yx_i) + G_p$ by definition and $G_p$ is $p$-divisible, it follows that $v_p(yx_i)$ is $p$-divisible in $G/G_p$ if and only if $\vn(yx_i)$ is $p$-divisible in $G$. It follows that the convex subgroup $\Delta$ of $G$ generated by $\vn(x)$ fulfils $|\Delta/p\Delta| \leq p^n$, since $\Delta/p\Delta \subseteq \{\vn(x_i) +p\Delta \mid i = 1,\ldots,p^n\}$. Hence $\Delta \subseteq G_{(p,n)}$ and in particular $\vn(x) \in G_{(p,n)}$, thus $x \in S_{(p,n)}$.\medskip
		
		Now we show for $x \in S_{p,n}$ that $K \models \psi_{(p,n)}(x)$. Note that $v_K(x) \in G_{(p,n)}$ since $v_p$ is the composition of $v_K$ with the projection from $G$ to $G/G_{(p,0)}$. We make the following case distinction:
		
		\begin{itemize}
			\item[(i)] If $K \models \varphi_p(x)$, then for all $y \in K$ fulfilling the part left of the implication in the first line of $\psi_{(p,n)}(x)$ it is $0 \leq v_p(y) \leq v_p(x)$. This implies $v_K(y) \in G_{(p,0)}$, $v_K(y) \in v_K(x) + G_{(p,0)}$  or $0 < v_K(y) < v_K(x)$. In total it follows that $v_K(y) \in G_{(p,n)}$.
			
			\item[(ii)] If $K \models \neg\varphi_p(x)$, then for all $y \in K$ fulfilling the part left of the implication in the second line of $\psi_{(p,n)}(x)$ we analogously obtain $v_K(y) \in G_{(p,0)}$, $v_K(y) \in v_K(x) + G_{(p,0)}$  or $v_K(x) < v_K(y) < 0$, so again $v_K(y) \in G_{(p,n)}$.
		\end{itemize}
		Since the right side of the implication in both lines is the same, it now suffices to show that for all $y \in K$ with $v_K(y) \in G_{(p,n)}$ there is an $i \in \{1,\ldots,p^n\}$ and a $z \in K$ such that $K \models \varphi_p(\frac{x_iy}{z^p}) \wedge \varphi_p(\frac{z^p}{x_iy})$. By choice of the parameters $x_i$ there is for every such $y$ an $\ell \in \{1,\ldots,p^n\}$ such that $-v_K(y) \in v_K(x_\ell) + p(G_{(p,n)})$. Now it is $v_K(x_\ell y) \in p(G_{(p,n)})$, hence we can choose $z \in K$ with $pv_K(z) = v_K(x_\ell y)$. Hence $v_K(\frac{x_iy}{z^p}) = v_K(\frac{z^p}{x_iy}) = 0$ and so is their projection to $G/G_{(p,0)}$, i.e.\ $v_p(\frac{x_iy}{z^p}) = v_p(\frac{z^p}{x_iy}) = 0$. We now proved that for $x \in S_{(p,n)}$ it follows that $K \models \psi_{(p,n)}(x)$.\\
		
		As a result the formula
		\[ \varphi_{(p,n)}(x) := \varphi_p(x) \vee \psi_{(p,n)}(x) \]
		defines $v_{(p,n)}$.\\

		\noindent(2): Assume for contradiction that $v_\Delta$ is a definable henselian valuation and $\Delta \not= G_{(p,n)}$ for all $(p,n)\in \{(p,i) \in \NN\times\NN_0 \mid p\textrm{ prime}, i \leq n_p\}$. Since by \cite[Corollary 3.25]{vandenDries14} the value group of a henselian valuation is stably embedded, $\Delta$ is a definable convex subgroup of $G$ and by \Cref{fact:dfblConvSubgp_pRegularCondition} there is a prime $p$ such that for all convex subgroups $\Delta_1, \Delta_2 \subseteq G$ with $\Delta_1 \subsetneq \Delta \subsetneq \Delta_2$ the quotient $\Delta_2/\Delta_1$ is not $p$-regular. Let $p_0$ be such a prime. As for $0 \leq i \leq n_{p_0}$ it is $\Delta \not= G_{(p_0,i)}$, we are in one of the following cases:
		\begin{itemize}
			\item[(i)] $\Delta \subsetneq G_{(p_0,0)} = G_{p_0}$: Choose $\Delta_2 = G_{(p_0,0)}, \Delta_1 = \{0\}$, then $\Delta_2/\Delta_1 = \Delta_2$ is $p_0$-divisible and hence in particular $p_0$-regular. Hence this case leads to a contradiction.
			
			\item[(ii)] $G_{(p_0,0)} \subsetneq \Delta \subsetneq G_{(p_0,n_{p_0})} = G$: Then the set $\{i \mid G_{(p_0,i)} \subsetneq \Delta\}$ has a maximum $m$. It follows that $G_{(p_0,m)} \subsetneq \Delta \subsetneq G_{(p_0,m+1)}$. Note that it is $|G_{(p_0,m+1)}/p_0G_{(p_0,m+1)}| = {p_0}^{m+1}$ since $G_{(p_0,m)} \subsetneq G_{(p_0,m+1)}$.
			
			We deduce that $G_{(p_0,m+1)}/G_{(p_0,m)}$ is $p_0$-regular:\newline			
			If $G_{(p_0,m+1)}/G_{(p_0,m)}$ is not $p_0$-regular, then there is a non-trivial convex subgroup $\Delta' \subsetneq G_{(p_0,m+1)}/G_{(p_0,m)}$ such that $(G_{(p_0,m+1)}/G_{(p_0,m)})/\Delta'$ is not $p_0$-divisible. But that means there exists $G_{(p_0,m)} \subsetneq \Delta'' \subsetneq G_{(p_0,m+1)}$ with $\Delta' = \Delta''/G_{(p_0,m)}$ and $G_{(p_0,m+1)}/\Delta''$ not $p_0$-divisible.
			
			We will now show that then $|\Delta''/p_0\Delta''| < {p_0}^{m+1}$:\newline
			It is $\Delta''/p_0G_{(p_0,m+1)} \subseteq G_{(p_0,m+1)}/p_0G_{(p_0,m+1)}$. Assume these set are equal, then for every $g \in G_{(p_0,m+1)}$ there is a $\delta \in \Delta''$ such that $g + p_0G_{(p_0,m+1)} = \delta + G_{(p_0,m+1)}$. Because $G_{(p_0,m+1)}/\Delta''$ was not $p_0$-divisible there is $g_0 \in G_{(p_0,m+1)}$ such that $g_0 + \Delta''$ is can not be divided by $p_0$ in $G_{(p_0,m+1)}/\Delta''$ $(*$).
			
			Choose $\delta_0 \in \Delta''$ with $\delta_0 + p_0G_{(p_0,m+1)} = g_0 + p_0G_{(p_0,m+1)}$. Then $g_0 - \delta_0 \in p_0G_{(p_0,m+1)}$, i.e.\ there is $g_1 \in G_{(p_0,m+1)}$ with $p_0g_1 = (g_0-\delta_0)$. But now
			\begin{align*}
				p_0 (g_1 + \Delta'') &= (p_0g_1) + \Delta''\\
				&= (g_0-\delta_0) + \Delta''\\
				&= g_0 +\Delta''.
			\end{align*}
			Contradiction to $(*)$, hence $\Delta''/p_0G_{(p_0,m+1)} \subsetneq G_{(p_0,m+1)}/p_0G_{(p_0,m+1)}$. In particular this means $|\Delta''/p_0G_{(p_0,m+1)}| < {p_0}^{m+1}$.
			
			Now for $\delta_1,\delta_2 \in \Delta''$ it is $\delta_1 -\delta_2 \in p_0G_{(p_0,m+1)}$ if and only if $\delta_1 - \delta_2 \in p_0\Delta''$ as $\Delta'' \cap p_0G_{(p_0,m+1)} = p_0\Delta''$. This follows because $\Delta''$ is a convex subgroup of $G_{(p_0,m+1)}$ and if for some $\delta \in \Delta'', g \in G_{(p_0,m+1)}$ it is $pg = \delta$, then $|g| < p_0|g| = |\delta|$, thus already $g \in \Delta''$. This yield a bijection from $\Delta''/p_0G_{(p_0,m+1)}$ to $\Delta''/p_0\Delta''$, hence $|\Delta''/p_0\Delta''| < {p_0}^{m+1}$.
			
			But $|\Delta''/p_0\Delta''|$ must be a power of $p_0$, hence it follows already that $|\Delta''/p_0\Delta''| \leq {p_0}^m$. Contradiction as $G_{(p_0,m)} \subsetneq \Delta''$ was the maximal convex subgroup of $G$ with $|G_{(p_0,m)}/p_0G_{(p_0,m)}| \leq {p_0}^{m}$.
			
			So  this case too leads to a contradiction with \Cref{fact:dfblConvSubgp_pRegularCondition}.
		\end{itemize}
		So the assumption of $v_\Delta$ being a definable henselian valuation with $\Delta \not= G_{(p,n)}$ for all $(p,n)\in \{(p,i) \in \NN\times\NN_0 \mid p\textrm{ prime}, i \leq n_p\}$ leads to a contradiction and must therefore be false. As a result all definable henselian valuations appear in the classification.
	\end{proof}

	\begin{remk}\label{rmk:notInjective}
		The map from \Cref{thm:ClassificationOfDefinableHenselianValuations} is in general not injective. On one hand, for different $p$ the $v_{(p,n)}$ can coincide through the same convex component adding non-$p$-divisible elements for multiple $p$, on the other hand for a single $p$ some $i \leq n_p$ might get skipped. Consider $\RR((\ZZ + \pi\ZZ))$ with the usual ordering (from the real numbers) of $\ZZ + \pi\ZZ$. Then $n_p = 2$ for all primes $p$, but there are no non-trivial proper convex subgroups, so $v_{(p,1)} = v_{(p,2)}$ for all primes $p$.
	\end{remk}

\section{Further work}

	In the proof of \Cref{thm:ClassificationOfDefinableHenselianValuations} the only property of the field $K$ that was necessary to show the definability of the valuations $v_{(p,n)}$ was almost real closedness. For any almost real closed field $K$ the valuations $v_{(p,n)}$ for a prime $p$ and $n \in \NN_0$ are therefore a subset of the definable henselian valuations on $K$. Is there a characterisation of all definable henselian valuations on an almost real closed field $K$ with strongly NIP value group $\vn(K^\times)$?
	

	\bigskip
	
	\noindent\textsc{Lothar Sebastian Krapp}:\smallskip\newline
	Institut für Vergleichende Sprachwissenschaft, Universit\"{a}t Z\"{u}rich, Z\"{u}rich, Switz-erland \& Fachbereich Mathematik und Statistik, Universit\"{a}t Konstanz, Konstanz, Germany\bigskip
	
	\noindent\textsc{Salma Kuhlmann}:\smallskip\newline
	Fachbereich Mathematik und Statistik, Universit\"{a}t Konstanz, Konstanz, Germany\bigskip
	
	\noindent\textsc{Lasse Vogel}:\smallskip\newline
	\texttt{lasse.vogel@uni-konstanz.de}\newline
	Fachbereich Mathematik und Statistik, Universit\"{a}t Konstanz, Konstanz, Germany
	

\begin{thebibliography}{99}
		\bibitem{DelonFarre96} \textsc{F. Delon} and \textsc{R. Farr\'{e}}, `Some model theory of almost real closed fields', \textsl{J. Symb. Log.} 61 (1996) 1121--1152.
		
		\bibitem{DittmannJahnkeKrappKuhlmann23} \textsc{P. Dittmann}, \textsc{F. Jahnke}, \textsc{L.S. Krapp} and \textsc{S. Kuhlmann}, `Definable valuations on ordered fields', \textsl{Model Theory} 2 (2023) 101--120.
		
		\bibitem{vandenDries14} \textsc{L. van den Dries}, `Lectures on the Model Theory of Valued Fields', \textsl{Model Theory in Algebra, Analysis and Arithmetic} (eds H. D. Macpherson and C. Toffalori), Lecture Notes in Mathematics 2111, (Springer, Berlin, 2014) 55--157.
		
		\bibitem{EnglerPrestel05} \textsc{A.J. Engler} and \textsc{A. Prestel}, \textsl{Valued Fields}, Springer Monogr. Math. (Springer, Berlin, 2005).
		
		\bibitem{FehmJahnke17} \textsc{A. Fehm} and \textsc{F. Jahnke}, `Recent progress on definability of Henselian valuations', Ordered Algebraic Structures and Related Topics, \textsl{Contemp. Math.} 697 (eds F. Broglia, F. Delon, M. Dickmann, D. Gondard-Cozette and V. A. Powers; Amer. Math. Soc., Providence, RI, 2017) 135--143.
		
		\bibitem{JahnkeSimonWalsberg17} \textsc{F. Jahnke}, \textsc{P. Simon} and \textsc{E. Walsberg}, `Dp-minimal valued fields', \textsl{J. Symb. Log} 82 (2017) 151--165.
		
		\bibitem{Johnson18} \textsc{W. Johnson}, `The canonical topology on dp-minimal fields', \textsl{J. Math. Log.} 18 (2018) Article ID: 1850007.
		
		\bibitem{Johnson23} \textsc{W. Johnson}, `The classification of dp-minimal and dp-small fields', \textsl{J. Eur. Math Soc. (JEMS)} 25 (2023) 467--513.
		
		\bibitem{KnebuschWright76} \textsc{M. Knebusch} and \textsc{M. J. Wright}, `Bewertungen mit reeler Henselisierung', \textsl{J. Reine Angew. Math.} 286-287 (1976) 314--321.
		
		\bibitem{KrappKuhlmannLehericy21} \textsc{L. S. Krapp}, \textsc{S. Kuhlmann} and \textsc{G. Leh\'{e}ricy}, `Ordered fields dense in their real closure and definable convex valuations', \textsl{Forum Math.} 33 (2021) 953--972.
		
		\bibitem{KrappKuhlmannLehericy21-2} \textsc{L. S. Krapp}, \textsc{S. Kuhlmann} and \textsc{G. Leh\'{e}ricy}, `Strongly NIP almost real closed fields', \textsl{Math. Log. Quart.} 67 (2021) 321--328.

		\bibitem{Shelah14} \textsc{S. Shelah}, `Strongly dependent theories', \textsl{Israel J. Math.} 204 (2014) 1--83.
		
		\bibitem{Simon15} \textsc{P. Simon}, \textsl{A Guide to NIP Theories}, Lecture Notes in Logic (Camebridge University Press, Camebridge, 2015).
	\end{thebibliography}
\end{document}